# Friable averages of oscillating arithmetic functions

R. de la Bretèche & G. Tenenbaum

*To the memory of Eduard Wirsing
whose profound insights will continue
fertilising our field.*

**Abstract.** We evaluate friable averages of arithmetic functions whose Dirichlet series is analytically close to some negative power of the Riemann zeta function. We obtain asymptotic expansions resembling those provided by the Selberg-Delange method in the non-friable case. An application is given to summing truncated versions of such functions.

## 1. Introduction and statements of results

Let $P^+(n)$ denote the largest prime factor of an integer $n$, with the convention that $P^+(1) = 1$. Given $y \geqslant 1$, an integer $n$ is said to be $y$-friable if $P^+(n) \leqslant y$. Let $S(x,y)$ stand for the set of $y$-friable integers not exceeding $x$ and let us write traditionally $\Psi(x,y) := |S(x,y)|$. Mean values of multiplicative functions over the set $S(x,y)$ attracted a lot of attention during the last decades—see, e.g., [3], [4], [6], [15], [16], [17]. However few articles deal with oscillating summands. Alladi [1], Hildebrand [7], [8], and Tenenbaum [12] consider the case of the Möbius function. Hildebrand's work [8] also deals with the function $e^{i\vartheta\Omega(n)}$, where $\vartheta \in \mathbb{R}$ and $\Omega(n)$ denotes the total number of prime factors, counted with multiplicity, of an integer $n$. In this work, we aim at handling a more general situation.

For given $\beta > 0$, $\mathfrak{c} \in ]0,1[$, $\delta > 0$, $\kappa > 0$, such that $\beta + \delta < 3/5$, we consider the classes $\mathcal{E}_\kappa^\pm(\beta, \mathfrak{c}, \delta)$ of those Dirichlet series $\mathcal{F}(s)$ convergent for $\sigma = \Re s > 1$ and which may be represented on this half-plane in the form

$$\mathcal{F}(s) = \zeta(s)^{\pm\kappa} \mathcal{B}(s) \tag{1.1}$$

where the series $\mathcal{B}(s) := \sum_{n \geqslant 1} b(n)/n^s$ may be holomorphically continued to a domain

$$\mathcal{D}(\beta, \mathfrak{c}, \delta) := \left\{ s \in \mathbb{C} : \sigma > 1 - \mathfrak{c}/\{1 + \log^+ |\tau|\}^{(1-\delta-\beta)/(\beta+\delta)} \right\} \tag{1.2}$$

where $\sigma := \Re e\, s$, $\tau := \Im s$, $\log^+ t = \max(\log t, 0)$ $(t > 0)$, and satisfies the conditions

$$\mathcal{B}(s) \ll \{1 + |\tau|\}^{1-\delta} \qquad (s \in \mathcal{D}(\beta, \mathfrak{c}, \delta)), \tag{1.3}$$

$$\mathcal{B}(s, y) := \sum_{P^+(n) \leqslant y} \frac{b(n)}{n^s} = \mathcal{B}(s) + O\left(\frac{1}{L_{\beta+\delta}(y)}\right) \quad \begin{pmatrix} y \geqslant 2,\ \sigma > 1 - \mathfrak{c}/(\log y)^{\beta+\delta}, \\ |\tau| \leqslant L_{\beta+\delta}(y) \end{pmatrix}, \tag{1.4}$$

with

$$L_r(y) := e^{(\log y)^r} \quad (r > 0,\ y \geqslant 2). \tag{1.5}$$

We then consider the class $\mathcal{H}(\kappa, \kappa_0; \beta, \mathfrak{c}, \delta)$ of those arithmetic functions $f$ whose associated Dirichlet series $\mathcal{F}(s)$ belongs to $\mathcal{E}_\kappa^-(\beta, \mathfrak{c}, \delta)$ and furthermore possesses a majorant series $\mathcal{F}^\dagger(s) = \zeta(s)^{\kappa_0} \mathcal{B}^\dagger(s)$ in $\mathcal{E}_{\kappa_0}^+(\beta, \mathfrak{c}, \delta)$. We note that the approach implemented in this work can readily be extended to the case where the powers of the Riemann zeta function appearing in (1·1) are replaced by products of powers of Dedekind zeta functions—see [6] for details. Moreover, it is plain that all our results below generalise *mutatis mutandis* to the case when $\kappa$ is an arbitrary complex number. However we chose to focus on the real case in order to avoid some technicalities.

We thus aim at establishing estimates for sums

$$M(x, y; f) := \sum_{n \in S(x,y)} f(n)$$

of real arithmetic functions $f$ with Dirichlet series

$$\mathcal{F}(s) := \sum_{n \geqslant 1} \frac{f(n)}{n^s} \qquad (\Re e\, s > 1)$$

belonging to $\mathcal{H}(\kappa, \kappa_0; \beta, \mathfrak{c}, \delta)$ for arbitrary values of the relevant parameters in the prescribed ranges.

The quantities arising in our formal statements below depend on the theory of delay differential equations. We describe the behaviour of the relevant functions in section 3.

Let $h_\kappa$ denote the unique continuous solution on $[0,\infty[$ of the equation

$$vh'_\kappa(v) = \kappa h_\kappa(v-1) \qquad (v>1) \tag{1.6}$$

with initial condition

$$h_\kappa(v) = 1 \qquad (0 \leqslant v \leqslant 1). \tag{1.7}$$

Define further

$$M(x;f) := M(x,x;f), \quad A(x,y;f) := x\int_{-\infty}^{\infty} h_\kappa(u-v)\,\mathrm{d}\Big(\frac{M(y^v;f)}{y^v}\Big). \tag{1.8}$$

Here and in the sequel, we write systematically $u := (\log x)/\log y$.

Our first result furnishes a range in which $A(x,y;f)$ is a good approximation to $M(x,y;f)$. It involves a rapidly decreasing function $R_\kappa$, precisely defined in §3.2 and satisfying (1.15) below.

**Theorem 1.1.** *Let* $\beta > 0$, $\mathfrak{c} > 0$, $\delta > 0$, $\kappa > 0$, $\kappa_0 > 0$, $\beta + \delta < 3/5$. *Then, for suitable* $\varepsilon = \varepsilon(\beta,\delta) > 0$ *and uniformly for* $f \in \mathcal{H}(\kappa,\kappa_0;\beta,\mathfrak{c},\delta)$ *and* $(x,y)$ *in the range*

$$(G_\beta) \qquad\qquad x \geqslant 3, \quad \exp\{(\log x)^{1-\beta}\} \leqslant y \leqslant x,$$

*we have*

$$M(x,y;f) = A(x,y;f) + O\Big(\frac{xR_\kappa(u)}{L_\varepsilon(y)}\Big). \tag{1.9}$$

Let $\nu := \lfloor\kappa\rfloor$, $\vartheta := \langle\kappa\rangle = \kappa - \nu$. We aim at some asymptotic expansion for $A(x,y;f)$ in the spirit of results provided by the Selberg-Delange method—see [14; ch. II.5 & II.6]. However this is by no means a straightforward consequence of (1.8). The finer behaviour will be described in terms of the function $\varphi_\kappa$, continuous solution on $\mathbb{R}^+$, differentiable on $[1,\infty[$, to the delay-differential equation

$$v\varphi'_\kappa(v) + \vartheta\varphi_\kappa(v) - \kappa\varphi_\kappa(v-1) = 0, \tag{1.10}$$

with initial condition

$$\varphi_\kappa(v) = \frac{v^{-\vartheta}}{\Gamma(1-\vartheta)} \quad (0 < v \leqslant 1). \tag{1.11}$$

It is worthwhile to observe right away that, if $\kappa = \nu \in \mathbb{N}^*$, then $\varphi_\kappa = h_\kappa$.

It follows from above that $\varphi_\kappa$ is $\mathcal{C}^\infty$ on $\mathbb{R}\smallsetminus\mathbb{N}$ and $\mathcal{C}^j$ on $]j,\infty[$ for all $j \geqslant 0$. When $\kappa \in \mathbb{N}^*$, the discontinuity of $\varphi_\kappa^{(j)}$ at $m \in [1,j]$ is of the first kind. We may hence continue $\varphi_\kappa^{(j)}$ by right-continuity at $m$. We write

$$\delta_{\kappa,m,j} = \varphi_\kappa^{(j)}(m) - \varphi_\kappa^{(j)}(m-0) \quad (\kappa \in \mathbb{N}^*,\, 1 \leqslant m \leqslant j). \tag{1.12}$$

It will be seen that

$$\psi_\kappa := \varphi_\kappa^{(\nu)} \tag{1.13}$$

occurs naturally in our study. A description of the asymptotic behaviour of $\psi_\kappa$ and its derivatives may be obtained from the general theory developed in [10]: we provide the necessary details in subsection 3.2.

For $v > 0$, define $\xi(v)$ as the unique non zero real solution of the equation $\mathrm{e}^\xi = 1+v\xi$, and put $\xi(1) := 0$. We have [14; lemma III.5.11][1]

$$\xi(v) = \log(v\log v) + O\Big(\frac{\log_2 v}{\log v}\Big) \qquad (v \geqslant 3),$$

and an asymptotic expansion may be derived through standard techniques. Define further $\zeta_0(v)$ as the solution of the equation $\mathrm{e}^\zeta = 1 - v\zeta$ with largest negative imaginary part. By [10; lemma 1], we have

$$\zeta_0(v) = \xi(v) + \frac{\pi^2}{2\xi(v)^2} - i\frac{\pi\xi(v)}{\xi(v)-1} + O\Big(\frac{1}{\xi(v)^3}\Big) \qquad (v \geqslant 2). \tag{1.14}$$

---

1. Here and throughout we denote by $\log_k$ the the $k$-fold iterated logarithm.

With this definition, we can state that the function $R_\kappa$ appearing in the statement of Theorem 1.1 and defined in §3.2 satisfies[2]

$$R_\kappa(v) \asymp \frac{1}{\sqrt{v}} \exp\left\{-\int_\kappa^v \zeta_0(t/\kappa)\,dt\right\} \qquad (v \geqslant 1). \tag{1.15}$$

Let $d\mu_{f,y}(t)$ be the measure on $\mathbb{R}$ with Laplace transform

$$\widehat{\mu_{f,y}}(s) := \int_\mathbb{R} e^{-ts}\,d\mu_{f,y}(t) = \frac{s^{1-\vartheta}\mathcal{F}(1+s/\log y)}{s+\log y} \qquad (\sigma > 0). \tag{1.16}$$

We then put, for $j \geqslant 0$, $y \geqslant 2$,

$$W_j(v,y;f) := \int_v^\infty (t-v)^j\,d\mu_{f,y}(t) \qquad (v \geqslant 0), \tag{1.17}$$

$$X_\ell(x,y;f) := \frac{(-1)^\ell}{(\ell-1)!}\int_0^{1/2} \varphi_\kappa^{(\ell)}(u-t)W_{\ell-1}(t,y;f)\,dt \quad (y^\ell < x \leqslant y^{\ell+1}). \tag{1.18}$$

It follows from estimates stated later—see in that order (5.1), (1.26), (5.15), and (5.10) *infra*—that

$$W_j(v,y;f) \ll \frac{1}{(\log y)^{j+\vartheta}L_{\beta+\delta/2}(y^v)} \qquad (v \geqslant 0,\, y \geqslant 2), \tag{1.19}$$

Moreover, as will be shown in §5.2, writing $\lambda := \log(x/y^\ell)$,

$$X_\ell(x,y;f) \ll \frac{\lambda^{1-\vartheta}}{(1+\lambda)(\log y)^\ell} + \frac{1}{(\log y)^{\ell+\vartheta}} \ll \frac{1}{(\log y)^\ell}. \tag{1.20}$$

Finally, for integer $J \geqslant 0$, real $y \geqslant 3$, put $\varepsilon_{J,y} := \{(2J+2)\log_2 y\}^{1/\beta}/\log y$ and consider the sets

$$\mathcal{D}_J(y) := \left\{u \geqslant 1 : \min_{1 \leqslant j \leqslant \min(u, J+1)} (u-j) \geqslant \varepsilon_{J,y}\right\}.$$

We can now state our second main result, where $\{a_j(f)\}_{j=0}^\infty$ is the sequence defined by the Taylor development

$$\frac{\mathcal{F}(s+1)}{s^\kappa(s+1)} = \sum_{j \geqslant 0} a_j(f)s^j \qquad (|s| < \mathfrak{c}).$$

**Theorem 1.2.** *Let*

$$\beta > 0,\ \mathfrak{c} > 0,\ \delta > 0,\ \kappa_0 > 0,\ \kappa > 0,\ \nu := \lfloor \kappa \rfloor,\ \beta + \delta < 3/5,\ J \in \mathbb{N}.$$

*Then, uniformly for $f \in \mathcal{H}(\kappa, \kappa_0; \beta, \mathfrak{c}, \delta)$, $(x,y) \in G_\beta$, $u = (\log x)/\log y \in \mathcal{D}_{J+\nu}(y)$, we have*

$$A(x,y;f) = x\sum_{0 \leqslant j \leqslant J} \frac{a_j(f)\psi_\kappa^{(j+1)}(u)}{(\log y)^{\kappa+j+1}} + O\left(\frac{xR_\kappa(u)(\log 2u)^{J+1}}{(\log y)^{\kappa+J+2}}\right). \tag{1.21}$$

*When $u \notin \mathcal{D}_{J+\nu}(y)$, $\ell < u \leqslant \ell+1$, and $\kappa \in \mathbb{N}^*$, the above formula persists if the quantity $xU_J(x,y;f)$ is added to the main term, with*

$$U_J(x,y;f) := \sum_{\ell \leqslant j \leqslant J+\nu+1} \frac{(-1)^{j+1}\delta_{\kappa,\ell,j}}{j!}W_j(u-\ell,y;f). \tag{1.22}$$

*If $u \notin \mathcal{D}_{J+\nu}(y)$, $\ell < u \leqslant \ell+1$, $\ell \leqslant J+\nu+1$, and $\kappa \in \mathbb{R}^* \smallsetminus \mathbb{N}^*$, formula (1.21) must be modified by restricting the summation to the (possibly empty) range $0 \leqslant j \leqslant \ell - \nu - 2$ and adding $xX_\ell(x,y;f)$ to the main term.*

We note that $\psi_1 = \omega$, Buchstab's function—see, e.g., [14; §III.6.2]—, so we recover [12; th. 2], with some further precision, in the special case $f = \mu$, the Möbius function.

---

2. Here and in the sequel we extend the meaning of the notation $f \asymp g$ to $|f| \asymp |g|$, hence relevant to complex quantities.

**Corollary 1.3.** *Let $\beta > 0$, $\mathfrak{c} > 0$, $\delta > 0$, $\kappa_0 > 0$, $\kappa > 0$, $\nu := \lfloor \kappa \rfloor$, $\beta + \delta < 3/5$. Then uniformly for $f \in \mathcal{H}(\kappa, \kappa_0; \beta, \mathfrak{c}, \delta)$, $(x, y) \in G_\beta$, $u \in \mathcal{D}_\nu(y)$, we have*

$$(1 \cdot 23) \qquad M(x, y; f) = \frac{x}{(\log y)^{\kappa+1}} \Big\{ \mathcal{B}(1) \psi'_\kappa(u) + O\Big(\frac{R_\kappa(u) \log(2u)}{\log y}\Big) \Big\}.$$

This describes the asymptotic behaviour unless $u \leqslant \nu + 1$ and $\langle u \rangle \leqslant \varepsilon_{\nu,y}$. When $\ell < u \leqslant \ell + 1$, $\ell \leqslant \nu + 1$, formula (1·23) must be modified according to the specifications described in the statement of Theorem 1.2. In particular, we have in all cases

$$(1 \cdot 24) \qquad M(x, y; f) \ll \frac{x}{(\log y)^{\min(\kappa+1, \ell)}} \qquad (y^\ell < x \leqslant y^{\ell+1}).$$

Note that, for $u = 1$, the Selberg-Delange method (see [14; ch. II.5]) furnishes

$$(1 \cdot 25) \quad M(x; f) = \frac{-\Gamma(\kappa+1) \sin(\pi \kappa) x}{\pi (1 + \log x)^{\kappa+1}} \Big\{ \mathcal{B}(1) + O\Big(\frac{1}{1 + \log x}\Big) \Big\} \quad (x \geqslant 1, \kappa \in \mathbb{R}^+ \smallsetminus \mathbb{N}),$$

$$(1 \cdot 26) \quad M(x; f) \ll \frac{x}{L_\beta(x)} \quad (\kappa \in \mathbb{N}^*).$$

Finally, we mention that the above results, the estimates of §§ 3.2-3.3, and of [15], [3], open the way to upper bounds that are uniform for $x \geqslant y \geqslant 2$. The following corollary is proved in § 7. We denote by $\mathcal{H}^*(\kappa, \kappa; \beta, \mathfrak{c}, \delta)$ the subclass of $\mathcal{H}(\kappa, \kappa; \beta, \mathfrak{c}, \delta)$ subject to the further condition that, for a suitable constant $C$, we have $\mathcal{B}(s, y) \ll \zeta(2\alpha_\kappa, y)^C$ and $1/\zeta(2\alpha_\kappa, y)^C \ll \mathcal{B}^\dagger(s, y) \ll \zeta(2\alpha_\kappa, y)^C$ when $\Re e\, s \geqslant \alpha_\kappa(x, y)$, defined as the saddle-point associated to $\zeta(s, y)^\kappa$.

**Corollary 1.4.** *Let $\beta > 0$, $\mathfrak{c} > 0$, $\delta > 0$, $\kappa_0 > 0$, $\kappa > 0$, $\beta + \delta < 3/5$, $1 < r < 3/2$. Then there exists a constant $c_0$ such that uniformly for $f \in \mathcal{H}^*(\kappa, \kappa; \beta, \mathfrak{c}, \delta)$, we have*

$$(1 \cdot 27) \qquad M(x, y; f) \ll M(x, y; f^\dagger) \Big\{ \frac{e^{-c_0 u/(\log 2u)^2}}{(\log y)^{\kappa+m-1}} + \frac{1}{L_r(y)} \Big\} \quad (x \geqslant y \geqslant 2),$$

where $m := \min(\lfloor u \rfloor, \kappa + 1)$.

## 2. Applications

### 2·1. *Weighted averages*

Consider the weighted analogue

$$(2 \cdot 1) \qquad m(x, y; f) := \sum_{n \in S(x, y)} \frac{f(n)}{n}$$

of $M(x, y; f)$. In some situations, it is convenient to have an estimate for (2·1) at our disposal, parallel to that following from Theorems 1.1 and 1.2. This is the purpose of the following statement in which $\mu^*_{f,y}$ denotes the measure with Laplace transform $\mathcal{F}(1 + s/\log y)/s^\vartheta$ and $\{a^*_j(f)\}_{j=0}^\infty$ is the sequence of Taylor coefficients of $\mathcal{F}(1+s)/s^\kappa$. For the sake of further reference we note that

$$(2 \cdot 2) \qquad a^*_j(f) = a_j(f) + a_{j-1}(f) \qquad (j \geqslant 0),$$

with the convention that $a_h(f) = 0$ if $h < 0$. We also define

$$(2 \cdot 3) \qquad W^*_j(v; f) := \int_v^\infty (t - v)^j \, \mathrm{d}\mu^*_{f,y}(t) \ll \frac{1}{L_{\beta+\delta/2}(y^v)(\log y)^{j+\vartheta}},$$

and, for $\ell < u \leqslant \ell + 1$,

$$(2 \cdot 4) \qquad X^*_\ell(x, y; f) := \frac{(-1)^\ell}{(\ell - 1)!} \int_0^{1/2} \varphi^{(\ell)}_\kappa(u - t) W^*_{\ell-1}(t, y; f) \, \mathrm{d}t.$$

Similarly to (1·20), writing $\lambda := \log(x/y^\ell)$, we have

(2·5) $$\frac{X_\ell^*(x,y;f) \ll \frac{\lambda^{1-\vartheta}}{(1+\lambda)(\log y)^\ell} + 1}{(\log y)^{\ell+\vartheta} \ll \frac{1}{(\log y)^\ell}}.$$

**Theorem 2.1.** *Let*

$$\beta > 0,\ \mathfrak{c} > 0,\ \delta > 0,\ \kappa_0 > 0,\ \kappa > 0,\ \nu := \lfloor \kappa \rfloor,\ \beta + \delta < 3/5,\ J \in \mathbb{N}.$$

*Then, uniformly for* $f \in \mathcal{H}(\kappa,\kappa_0;\beta,\mathfrak{c},\delta)$, $(x,y) \in G_\beta$, $u \in \mathcal{D}_{J+\nu}(y)$, *we have*

(2·6) $$m(x,y;f) = \sum_{0 \leqslant j \leqslant J} \frac{a_j^*(f)\psi_\kappa^{(j)}(u)}{(\log y)^{\kappa+j}} + O\left(\frac{R_\kappa(u)(\log 2u)^J}{(\log y)^{J+\kappa+1}}\right),$$

*When* $u \notin \mathcal{D}_{J+\nu}(y)$, $\ell < u \leqslant \ell+1$, *and* $\kappa \in \mathbb{N}^*$, *the above formula persists provided the quantity* $U_J^*(x,y;f)$ *is added to the main term, with*

(2·7) $$U_J^*(x,y;f) := \sum_{\ell \leqslant j \leqslant J+\nu} \frac{(-1)^{j+1}\delta_{\kappa,\ell,j}}{j!}W_j^*(u-\ell;f),$$

*where* $\delta_{\kappa,m,j}$ *is defined in* (1·12).

*If* $u \notin \mathcal{D}_{J+\nu}(y)$, $\ell < u \leqslant \ell+1$, $\ell \leqslant J+\nu+1$, *and* $\kappa \in \mathbb{R}^* \smallsetminus \mathbb{N}^*$, *formula* (2·6) *must be modified by restricting the summation to the (possibly empty) range* $0 \leqslant j \leqslant \ell - \nu - 2$ *and adding* $X_\ell^*(x,y;f)$ *to the main term.*

At the cost of a weakening of the error term, one can take advantage of the rapid decrease of the density of friable integers as $u$ gets large in order to derive estimates valid without any restriction. For instance in the case of the Möbius function $\mu$, we have, uniformly for $x \geqslant y \geqslant 2$,

$$\sum_{n \in S(x,y)} \frac{\mu(n)}{n} = \frac{\omega(u)}{\log y}\int_1^{x/y} \frac{m(t)}{t}\,\mathrm{d}t + O\left(\frac{1}{(\log y)^2}\right),$$

where $\omega$ denotes Buchstab's function and $m(t) := \sum_{n \leqslant t} \mu(n)/n$.

### 2·2. *Truncated multiplicative functions*

In a recent preprint [2], Alladi & Goswami gave estimates for the summatory function

$$\sum_{n \leqslant x}(-k)^{\omega(n,y)}$$

where $k \in \mathbb{N}^*$, and $\omega(n,y) := \sum_{p\,|\,n,\,p \leqslant y} 1$. The results presented in section 1 enable to consider more generally $M(x;f_y)$ for $f \in \mathcal{H}(\kappa,\kappa_0;\beta,\mathfrak{c},\delta)$ with $\kappa > 0$ and $f_y$ is defined as $f_y(n) = f(m)$ if $m$ is the largest $y$-friable divisor of $n$. Thus $f_y$ is obtained from $f$ by truncating its values at all large primes, similarly to the process classically used for additive functions.

Some further notation is necessary to state our result.

For integer $J \geqslant 0$ and constant $b > 0$, we consider the sets

$$\mathcal{D}_J(b,y) := \left\{u \geqslant 1 : \min_{1 \leqslant j < \min(u,J+1)}(u-j) > 1/(\log y)^b\right\}.$$

Our estimate for $M(x;f_y)$ is stated below for $u \in \mathcal{D}_J(b,y)$ and suitable $b > 0$. This produces an estimate in which the discontinuities of $\psi_{\kappa+1}$ have no influence. Taking into account the contributions of the discontinuities described in Theorems 1.2 and 2.1 leads to a more complicated statement valid without restriction. We leave this to the reader.

**Theorem 2.2.** *Let*

$$\beta > 0,\ \mathfrak{c} > 0,\ \delta > 0,\ \kappa_0 > 0,\ \kappa > 0,\ \nu := \lfloor \kappa \rfloor,\ \beta + \delta < \tfrac{1}{2},\ J \in \mathbb{N}^*,\ b := (1-2\beta)/(1-\beta).$$

Then, uniformly for $f \in \mathcal{H}(\kappa, \kappa_0; \beta, \mathfrak{c}, \delta)$, $(x, y) \in G_\beta$, and $u \in \mathcal{D}_{J+\nu+1}(b, y)$, we have

$$(2\cdot 8) \qquad M(x; f_y) = x \sum_{0 \leqslant j \leqslant J} \frac{a_j(f)\psi_{\kappa+1}^{(j)}(u)}{(\log y)^{\kappa+j+1}} + O\bigg(\frac{xR_\kappa(u)(\log 2u)^{J+1}}{(\log y)^{J+\kappa+2}}\bigg).$$

In the special case $f(n) := (-k)^{\omega(n)}$, where $\omega(n)$ now counts without multiplicity the total number of prime factors of $n$, formula (2·8) represents a significant improvement over the corresponding estimate in [2]: it is valid in a much larger range since [2; th. 5.2] requires $\beta$ to be taken arbitrarily small, and it also furnishes an expansion according to negative powers of $\log y$ whereas only the dominant term is provided in [2]. Note that for this particular function $a_0(f)$ vanishes whenever $k = p + 1$ for some prime $p$. This accounts for the dichotomy put forward in [2].

It is also noteworthy to remark that, under the assumptions of Theorem 2.2 and $u \gg \log_2 y$, an obvious modification of the proof of [5; th. 02] provides

$$M(x; f_y) = \frac{a_0(f)e^{-\gamma(\kappa+1)}x}{(\log y)^{\kappa+1}} + O\bigg(\frac{x}{(\log y)^{\kappa+2}}\bigg).$$

The main term agrees with that of (2·8) in view of (3·6) below. Hence (2·8) extends the scope of the above result by providing an expansion of the remainder term in a fairly large domain.

## 3. Solutions to delay-differential equations

### 3·1. The function $h_\kappa$

We continue $h_\kappa$ on $\mathbb{R}$ by setting $h_\kappa(v) = 0$ for $v < 0$, so that (1·6) still holds for $v \in \mathbb{R} \smallsetminus \{0, 1\}$. From the general theory displayed in [10], we know—see [10; theorem 1]—that

$$h_\kappa(v) \ll v^\kappa \qquad (v > 1),$$

so the Laplace transform

$$\widehat{h_\kappa}(s) := \int_0^\infty e^{-vs} h_\kappa(v)\,dv$$

converges for $\sigma := \Re e\, s > 0$. Arguing as in [14; § III.6.3], which corresponds to a similar situation in the case $\kappa = 1$, it is easy to show that $\big(s\widehat{h_\kappa}(s)\big)' = -\kappa e^{-s}\widehat{h_\kappa}(s)/s$, and then

$$(3\cdot 1) \qquad \widehat{h_\kappa}(s) = \frac{1}{s^{\kappa+1}\widehat{\varrho}(s)^\kappa} \qquad (\sigma > 0)$$

where $\varrho$ denotes Dickman's function, the solution of (1·6) when $\kappa$ is replaced by $-1$. We recall (see, e.g. [14; th. III.5.10]) that, writing $I(s) := \int_0^s (e^t - 1)\,dt/t$, we have

$$(3\cdot 2) \qquad \widehat{\varrho}(s) = e^{\gamma + I(-s)} \qquad (s \in \mathbb{C}).$$

By [10; eq. (3.4″)], we have

$$h_\kappa(v) = AF(v; 0, -\kappa) + 2\Re e \sum_{0 \leqslant n < N} A_n F_n(v; 0, -\kappa) + O(E_N(v))$$

$$E_N(v) \ll F_0(v; 0, -\kappa) \exp\bigg\{-\frac{2N(N+1)\pi^2 v}{\xi(v)^2}\bigg(1 + O\bigg(\frac{1}{\xi(v)}\bigg)\bigg)\bigg\} \qquad (v \geqslant 2)$$

where, here and in the sequel, we let $F(v; a, b)$ denote the exceptional solution, and $F_n(v; a, b)$ the $n$th fundamental solution,[3] of the delay differential equation

$$(3\cdot 3) \qquad vf'(v) + af(v) + bf(v-1) = 0.$$

Moreover,

$$A := h_\kappa(1)G(1; 0, -\kappa) + \kappa \int_0^1 h_\kappa(v) G(v+1; 0, -\kappa)\,dv = \lim_{v \to 0+} vG(v; 0, -\kappa),$$

$$A_n := h_\kappa(1)G_n(1; 0, -\kappa) + \kappa \int_0^1 h_\kappa(v) G_n(v+1; 0, -\kappa)\,dv$$

$$= \lim_{v \to 0+} vG_n(v; 0, -\kappa) \qquad (n \geqslant 0).$$

where $G(v; a, b)$ is the exceptional solution, and $G_n(v; a, b)$ is the $n$th fundamental solution, of the adjoint equation $vg'(v) + (1-a)g(v) - bg(v+1) = 0$.

---

3. These functions are precisely defined as convergent contour integrals in [10].

From the expression of $G(v; 0, -\kappa)$ given in [10; (2.7)], it can be deduced that

$$A = e^{-\gamma\kappa}/\Gamma(\kappa+1),$$

where $\gamma$ is Euler's constant and $\Gamma$ denotes Euler's function—see [13; (2.13)] for a similar computation relevant to a solution of (3·3) with $(a,b) = (1,-1)$. Since it follows from [10; th. 1] that $F(v) = v^\kappa + O(v^{\kappa-1})$ as $v \to \infty$, while $F(v)$ is indeed a polynomial of degree $\kappa$ whenever $\kappa$ is a positive integer, we derive, for any fixed integer $J \geqslant 0$,

$$(3\cdot 4) \qquad h_\kappa(v) = \sum_{0 \leqslant j \leqslant J} \frac{c_j v^{\kappa-j}}{\Gamma(\kappa+1-j)} + O\bigl(v^{\kappa-J-1}\bigr) \qquad (v \to \infty),$$

where $\{c_j\}_{j=0}^\infty$ is the sequence of Taylor coefficients at the origin of $e^{-\gamma\kappa - \kappa I(-s)} = \widehat{\varrho}(s)^{-\kappa}$.

### 3·2. *The functions $\varphi_\kappa$ and $\psi_\kappa$*

Recall (1·10) and (1·11). The Laplace transform

$$\widehat{\varphi_\kappa}(s) := \int_0^\infty \varphi_\kappa(v) e^{-sv} \, dv \qquad (\sigma > 0)$$

can be computed classically by showing from the delay-differential equation (1·10) that it satisfies a linear differential equation. We omit the details, which are similar to the computation of $\widehat{\varrho}(s)$—see, e.g., [14; th. III.5.10]. We obtain

$$(3\cdot 5) \qquad \widehat{\varphi_\kappa}(s) = \frac{1}{s^{\nu+1} \widehat{\varrho}(s)^\kappa} \qquad (\sigma > 0),$$

and observe that the inverse Laplace integral

$$\varphi_\kappa(v) = \frac{1}{2\pi i} \int_{1+i\mathbb{R}} \frac{e^{vs}}{s^{\nu+1} \widehat{\varrho}(s)^\kappa} \, ds \qquad (v \in \mathbb{R})$$

converges for all $v \neq 0$ since $\widehat{\varrho}(s) = 1/s + O(1/s^2)$ for $\Re e\, s = 1$. Moreover, we recover the fact that $\varphi_\kappa(v) = 0$ for $v < 0$ by moving the integration line to the right at infinity.

From (1·10), we see that $\psi_\kappa$ satisfies (3·3) with $(a,b) = (\kappa, -\kappa)$. By [10; th. 1], we have $F(v; \kappa, -\kappa) = 1$ for the exceptional solution of this delay-differential equation. Therefore, letting $\delta_{0j}$ denote Kronecker's symbol, we get for $j \geqslant 0$ that

$$(3\cdot 6) \qquad \psi_\kappa^{(j)}(v) = \delta_{0j} e^{-\gamma\kappa} + 2\Re e \sum_{0 \leqslant n < N} \lambda_{j,n} F_n(v; \kappa+j, -\kappa) + O(E_N(v)),$$

$$(3\cdot 7) \qquad E_N(v) \ll F_0(v; \kappa+j, -\kappa) \exp\left\{-\frac{2N(N+1)\pi^2 v}{\xi(v)^2}\left(1 + O\Bigl(\frac{1}{\xi(v)}\Bigr)\right)\right\},$$

$$(3\cdot 8) \qquad \psi_\kappa^{(j+1)}(v) \ll (\log v)^j R_\kappa(v) \qquad (v \to \infty),$$

for $v \geqslant 2$, $N \geqslant 1$, with $\lambda_{j,0} \neq 0$. The constant term in (3·6) has been identified by comparing saddle-point estimates of $F(v; -\kappa, \kappa)$ and $\psi_\kappa$. Similarly, using [10; th. 2, (2.17), (3.13)], it can be shown that $\lambda_{0,0} = \lambda_{1,0} = e^{-\gamma\kappa}$.

The behaviour of fundamental solutions $F(v; a, b)$ and $F_n(v; a, b)$ ($n \geqslant 0$) is closely linked to the complex roots of the equation $e^\zeta = 1 - v\zeta$. These form an infinite sequence located in the neighbourhood of the vertical line $\Re e\, \zeta = \xi(v)$. However, the equation fails to have a non-trivial real root. The two conjugate roots closer to the real axis have hence to be used in order to describe the behaviour of solutions to (3·3). This has been achieved, in a general context in [10]. With the notation of this last work, we denote these two roots by $\zeta_0(v)$ and $\zeta_{-1}(v) = \overline{\zeta_0(v)}$. We write $\xi_0(v) := \Re e\, \zeta_0(v)$ $(v \geqslant 1)$ and recall the asymptotic formula (1·14).

Define

$$(3\cdot 9) \qquad \Phi_\kappa(v,s) := \frac{\exp\{-vs - \kappa I(s)\}}{\sqrt{2\pi v(1 - 1/s)}} = \frac{e^{\gamma\kappa - vs}}{\widehat{\varrho}(-s)^\kappa \sqrt{2\pi v(1 - 1/s)}}.$$

The asymptotic behaviour of the fundamental solutions $F_n(v;a,b)$ has been precisely described in [10]. By [10; (2.1),(2.17)], we have

$$(3\cdot10) \quad F_0(v;\kappa+1,-\kappa) = \left\{1+O\left(\frac{1}{v}\right)\right\}\Phi_\kappa(v,\zeta_0(v/\kappa)),$$

$$(3\cdot11) \quad R_\kappa(v) := |F_0(v;\kappa+1,-\kappa)| \asymp \frac{e^{-v\xi_0(v/\kappa)}}{|\widehat{\varrho}(-\zeta_0(v/\kappa))|^\kappa \sqrt{v}} \asymp \frac{1}{\sqrt{v}}\exp\left\{-\int_\kappa^v \zeta_0(t/\kappa)\,dt\right\},$$

$$(3\cdot12) \quad \Re e\, F_0(v;\kappa+1,-\kappa) = R_\kappa(v)\left\{\cos\vartheta_\kappa(v)+O\left(\frac{1}{v}\right)\right\}$$

where $\vartheta_\kappa(v) := \arg\Phi_\kappa(v,\zeta_0(v/\kappa))$. Furthermore, it may be derived from [13; lemma 4], corresponding to the case $\kappa = 1$, that

$$\vartheta'_\kappa(v) = \frac{\pi\xi(v/\kappa)}{\xi(v/\kappa)-1} + O\left(\frac{1}{\xi(v)^3}\right) \quad (v \geqslant 2\kappa).$$

### 3·3. Link with positive convolution powers of the Dickman function

We end this section by mentioning that the functions $\Phi_\kappa$ and $\xi$ are also relevant for describing the behaviour of the function $\varrho_\kappa$ defined as the order $\kappa$ convolution power of $\varrho = \varrho_1$. Observing that $\varrho_\kappa$ satisfies the equation $v\varrho'_\kappa(v)-(\kappa-1)\varrho_\kappa(v)+\kappa\varrho_\kappa(v-1)=0$ ($v>1$) with initial condition $\varrho_\kappa(v) = v^{\kappa-1}/\Gamma(\kappa)$ for $0 < v \leqslant 1$, we also infer from [10] (see also [11]) that, for large $v$,

$$\varrho_\kappa(v) = \left\{e^{\gamma\kappa}+O\left(\frac{1}{v}\right)\right\}\Phi_\kappa(v,\xi(v/\kappa))$$
$$= \left\{1+O\left(\frac{1}{v}\right)\right\}\sqrt{\frac{\xi'(v/\kappa)}{2\kappa\pi}}\exp\left\{\gamma\kappa - \int_\kappa^v \xi(t/\kappa)\,dt\right\},$$

and so, in view of (1·14), (3·9) and (3·10), writing $H(v) := \exp\{v/(\log 2v)^2\}$,

$$(3\cdot13) \quad R_\kappa(v) = \varrho_\kappa(v)\exp\left\{\frac{-\pi^2 v}{2\xi(v)^2} + O\left(\frac{v}{\xi(v)^3}\right)\right\} = \frac{\varrho_\kappa(v)}{H(v)^{\pi^2/2+o(1)}}.$$

Recall that, in the case $\mathcal{F}(s) = \zeta(s)^\kappa \mathcal{B}(s) \in \mathcal{E}_\kappa^+(\beta,\mathfrak{c},\delta)$, we have, by [6; th. 1.1],

$$(3\cdot14) \quad M(x,y;f) = x\varrho_\kappa(u)(\log y)^{\kappa-1}\left\{\mathcal{B}(1) + O\left(\frac{\log 2u}{\log y} + \frac{1}{(\log y)^\kappa}\right)\right\},$$

uniformly for $(x,y) \in G_\beta$, with $u := (\log x)/\log y$.

## 4. Proof of Theorem 1.1

### 4·1. Auxiliary estimates

**Lemma 4.1 ([13]).** For $v \geqslant 1$, $\Re e\, s = -\xi_0(v)$, we have

$$(4\cdot1) \quad \left|\zeta_0(v)\widehat{\varrho}(-\zeta_0(v))\right| \ll |s\widehat{\varrho}(s)|.$$

*Proof.* The bound (4·1) coincides with [13; lemma 8]. □

Our second lemma states an estimate which may be proved similarly to [13; lemma 10], the details being left to the reader. We use the notation

$$(4\cdot2) \quad \alpha_\kappa = \alpha_\kappa(x,y) := 1 - \xi_0(u/\kappa)/\log y \quad (x \geqslant y \geqslant 2),$$

$$(4\cdot3) \quad \zeta(s,y) := \prod_{p \leqslant y}(1-1/p^s)^{-1} \quad (\Re e\, s > 0, y \geqslant 2),$$

and, for $b > 0$, define the domain

$$(H_b) \quad x \geqslant 3, \ \exp\{(\log_2 x)^b\} \leqslant y \leqslant x.$$

**Lemma 4.2.** *Let $b > 5/3$. Uniformly for $(x, y) \in (H_b)$, we have*

$$(4 \cdot 4) \qquad x^{\alpha_\kappa} \zeta(\alpha_\kappa, y)^\kappa \ll x \varrho_\kappa(u) e^{O(u/\xi(u)^4)} (\log y)^\kappa.$$

We next restate [14; lemma III.5.16] in the following weaker form.

**Lemma 4.3.** *Let $3/5 < b < 3/2$. Then, for $0 < \varepsilon < 1/2 - b/3$, we have*

$$(4 \cdot 5) \qquad \zeta(s, y) = \zeta(s) s_y \widehat{\varrho}(s_y) \Big\{ 1 + O\Big(\frac{1}{L_\varepsilon(y)}\Big) \Big\}$$

*with $s_y := (s - 1) \log y$ and uniformly in the range*

$$y \geqslant 2, \quad \sigma \geqslant 1 - 1/(\log y)^{2b/3 + 3\varepsilon/2}, \quad |\tau| \leqslant L_b(y).$$

Finally, we need an estimate for short sums of the coefficients of series in $\mathcal{E}_\kappa^+(\beta, \mathfrak{c}, \delta)$.

**Lemma 4.4.** *Let $\beta > 0$, $\mathfrak{c} > 0$, $\delta > 0$, $\beta + \delta < 3/5$, $\kappa_0 > 0$, and let $f^\dagger$ denote a non-negative arithmetic function with Dirichlet series in $\mathcal{E}_{\kappa_0}^+(\beta, \mathfrak{c}, \delta)$. Then the estimate*

$$(4 \cdot 6) \qquad \sum_{x < n \leqslant x + z} f^\dagger(n) \ll z (\log x)^{\kappa_0 - 1} + \frac{x}{L_{\beta + 2\delta/3}(x)}$$

*holds uniformly for $x \geqslant 3$, $0 \leqslant z \leqslant x$.*

*Proof.* This is a straight-forward consequence of the estimate for the summatory function of $f^\dagger$ achieved via the Selberg–Delange method, as displayed in [14; ch. II.5]. We note for further reference that (4·6) implies

$$(4 \cdot 7) \qquad f^\dagger(n) \ll n / L_{\beta + \delta/2}(n)^{10} \qquad (n \geqslant 1). \qquad \square$$

### 4·2. *Completion of the proof*

Let $\delta > 0$ be small, put $T := u^{2u} L_{\beta + \delta/2}(y)$, and let $\{f^\dagger(n)\}_{n=1}^\infty$ denote the sequence of the coefficients of the majorant series $\mathcal{F}^\dagger(s) \in \mathcal{E}_\kappa^+(\beta, \mathfrak{c}, \delta)$. We first apply Perron's formula (see, e.g., [14; th. II.2.3]) to get, for $(x, y) \in G_\beta$,

$$(4 \cdot 8) \qquad M(x, y; f) = \frac{1}{2\pi i} \int_{\alpha_\kappa - iT^2}^{\alpha_\kappa + iT^2} \frac{\mathcal{B}(s, y) x^s}{s \zeta(s, y)^\kappa} \, ds + \mathfrak{R},$$

with

$$\mathfrak{R} \ll \sum_{P(n) \leqslant y} \frac{x^{\alpha_\kappa} f^\dagger(n)}{n^{\alpha_\kappa} (1 + T^2 |\log(x/n)|)} \ll \frac{x^{\alpha_\kappa} \zeta(\alpha_\kappa, y)^\kappa}{T} + \sum_{|n - x| \leqslant x/T} f^\dagger(n) \ll \frac{x R_\kappa(u)}{L_{\beta + \delta/2}(y)},$$

by (4·4), (3·13), and (4·6).

Next, we check that, for sufficiently small $\delta$, we have $T \leqslant L_r(y)$ for some $r < 3/2$, and so apply formula (3·1) and Lemma 4.3 to derive, for $(x, y) \in (G_\beta)$, $|\tau| \leqslant T^2$, $s_y := (s - 1) \log y$, $L = L_\varepsilon(y)$,

$$(4 \cdot 9) \qquad \frac{x^s}{s \zeta(s, y)^\kappa} = \frac{x e^{u s_y} s_y \widehat{h_\kappa}(s_y)}{s \zeta(s)^\kappa} \Big\{ 1 + O\Big(\frac{1}{L^2}\Big) \Big\} = \frac{x e^{u s_y} \{1 + O(1/L^2)\}}{s \{\zeta(s) s_y\}^\kappa \widehat{\varrho}(s_y)^\kappa}.$$

We note that, when $\kappa \notin \mathbb{N}^*$, the singularity arising from the pole of $\zeta(s)$ at $s = 1$ is compensated by the corresponding zero of $s_y$, so that the main term is analytic on a zero-free region of the zeta function.

Applying (4·1) and (3·11), we see that the contribution of the last error term to the integral of (4·8) is

$$\ll \frac{x^{\alpha_\kappa}}{L^2} \int_0^{T^2} \frac{\{\log(2 + \tau)\}^\kappa \, d\tau}{|\zeta_0(u/\kappa) \widehat{\varrho}(-\zeta_0(u/\kappa))|^\kappa (1 + \tau)} \ll \frac{x R_\kappa(u) \sqrt{u} \{\log T\}^{\kappa + 1}}{L^2} \ll \frac{x R_\kappa(u)}{L}.$$

This is compatible with (1·9).

Taking (1·4) into account, the same argument yields

$$(4\cdot10) \quad M(x,y;f) = \frac{1}{2\pi i}\int_{\alpha_\kappa-iT^2}^{\alpha_\kappa+iT^2} \frac{\mathcal{B}(s)x^s}{s\zeta(s)^\kappa s_y^\kappa \widehat{\varrho}(s_y)^\kappa}\,\mathrm{d}s + O\Big(\frac{xR_\kappa(u)}{L}\Big) \quad (x \geqslant y \geqslant 2).$$

Considering (3·1) and since

$$(4\cdot11) \quad \frac{(s-1)\mathcal{B}(s)}{s\zeta(s)^\kappa} = \int_0^\infty \mathrm{e}^{-v(s-1)}\,\mathrm{d}\Big(\frac{M(\mathrm{e}^v;f)}{\mathrm{e}^v}\Big) \quad (\sigma > 1),$$

we deduce from the convolution theorem that $s_y\mathcal{B}(s)\widehat{h_\kappa}(s_y)/\{s\zeta(s)^\kappa\}$ is the Laplace transform of

$$J_y(t) := \mathrm{e}^t \int_0^\infty h_\kappa\Big(\frac{t}{\log y} - v\Big)\,\mathrm{d}\Big(\frac{M(y^v;f)}{y^v}\Big).$$

We plainly have $J_y(\log x) = A(x,y;f)$. Since $\widehat{J_y}(s)$ is holomorphic in any zero-free region of the zeta function, we may write

$$(4\cdot12) \quad J_y(\log x) = \frac{1}{2\pi i}\int_\mathcal{L} \widehat{J_y}(s)x^s\,\mathrm{d}s$$

where $\mathcal{L}$ is the broken line $[b-i\infty, b-iT^2, \alpha_\kappa-iT^2, \alpha_\kappa+iT^2, b+iT^2, b+i\infty]$ with $b := 1+1/\log x$. In order to prove (1·9), it therefore remains to show that the contribution of the complement of the segment $[\alpha_\kappa - iT^2, \alpha_\kappa + iT^2]$ is negligible.

On the horizontal segments, we have $s_y\widehat{h_\kappa}(s_y) \ll 1$ and, by standard bounds for the zeta function in the Vinogradov-Korobov region,

$$\mathcal{F}(s) = \mathcal{B}(s)/\zeta(s)^\kappa \ll \big(\log(|\tau|+2)\big)^\kappa.$$

The corresponding contribution is therefore $\ll x(\log T)^\kappa/T^2$, which is plainly acceptable.

Now, by [14; lemma III.5.12], we have

$$s_y\widehat{h_\kappa}(s_y) = 1 + O\Big(\frac{1+u\xi(u)}{s_y}\Big) \quad (|\tau| > T^2).$$

Thus the contribution of the vertical half-lines may be estimated by Perron's formula for the main term $x^s/\{s\zeta(s)^\kappa\}$ and by a direct bound for the remainder $\ll x\{1+u\xi(u)\}/(\tau^2\log y)$.

This completes the proof of (1·9).

## 5. Proof of Theorem 1.2

### 5·1. The case $\kappa \in \mathbb{N}^*$

When $\kappa = \nu \in \mathbb{N}^*$, we have, by (1·16) and (4·11),

$$(5\cdot1) \quad \mu_{f,y}(v) = M(y^v;f)/y^v,$$

$$(5\cdot2) \quad \frac{a_{j-\kappa-1}(f)}{(\log y)^j} = \frac{(-1)^j}{j!}\int_0^\infty v^j\,\mathrm{d}\mu_{f,y}(v) \quad (j \geqslant 0),$$

with the convention that $a_h(f) = 0$ if $h < 0$.

We may discard the case $u = 1$ since the required formula is then a straightforward consequence of the Selberg–Delange method, as stated in [14; th. II.5.2].

When $u \in \mathcal{D}_{J+\nu}(y)$, $u > 1$, we apply (1·9) and (5·1) and write

$$\int_0^u \varphi_\kappa(u-v)\,\mathrm{d}\mu_{f,y}(v) = I_1 + I_2 + O\Big(\frac{1}{L_{\beta+\delta/2}(x)}\Big),$$

where the $I_j$ correspond respectively to the integration ranges $[0,\tfrac{1}{2}\varepsilon_y]$, $]\tfrac{1}{2}\varepsilon_y, u-\tfrac{1}{2}]$; here and throughout, we write for simplicity $\varepsilon_y := \varepsilon_{J+\nu,y}$.

To evaluate $I_1$ we may use the fact that $\varphi_\kappa$ belongs the class $\mathcal{C}^{J+\kappa+2}$ on $[u - \tfrac{1}{2}\varepsilon_y, u]$. For $0 \leqslant v \leqslant \tfrac{1}{2}\varepsilon_y$ we hence have the Taylor-Lagrange expansion

$$\varphi_\kappa(u-v) = \sum_{0 \leqslant j \leqslant J+\kappa+1} \frac{(-v)^j \varphi_\kappa^{(j)}(u)}{j!} + \mathfrak{R}_0, \tag{5.3}$$

with

$$\mathfrak{R}_0 := \frac{(-1)^{J+\kappa+2}}{(J+\kappa+1)!} \int_0^v (v-t)^{J+\kappa+1} \psi_\kappa^{(J+2)}(u-t)\,\mathrm{d}t.$$

Therefore

$$I_1 = \sum_{0 \leqslant j \leqslant J+\kappa+1} \frac{(-1)^j \varphi_\kappa^{(j)}(u)}{j!(\log y)^{\kappa+j+1}} \int_0^{\varepsilon_y/2} v^j \,\mathrm{d}\mu_{f,y}(v) + \mathfrak{S}_0, \tag{5.4}$$

with

$$\mathfrak{S}_0 := \frac{(-1)^{J+\kappa+2}}{(J+\kappa+1)!} \int_0^{\varepsilon_y/2} \psi_\kappa^{(J+2)}(u-t) \int_t^{\varepsilon_y/2} (v-t)^{J+\kappa+1}\,\mathrm{d}\mu_{f,y}(v)\,\mathrm{d}t.$$

We first observe that, by (1.26), we may extend the integrals in (5.4) to infinity involving a remainder absorbable by that of (1.21). To this extent and in view of (5.2), the sum may be replaced by

$$\sum_{0 \leqslant j \leqslant J} \frac{a_j(f)\psi_\kappa^{(j)}(u)}{(\log y)^{\kappa+j+1}}.$$

Next, noting that, for $j > \nu$, $u \geqslant 1, 0 \leqslant v \leqslant u - \tfrac{1}{2}$,

$$\psi_\kappa^{(j)}(u-v) \ll R_\kappa(u-v)(\log 2u)^{j-1} \ll R_\kappa(u)(\log 2u)^{j-1} \mathrm{e}^{v\xi_0(u/\kappa)}, \tag{5.5}$$

which readily follows, as in [14; cor. III.5.15], from (3.8) and (3.11), we see that $\mathfrak{S}_0$ may be absorbed by the remainder of (1.21).

The integral $I_2$ can be handled as an error term. Indeed, it suffices to perform partial summation and apply (1.26) to see that $I_2$ does not exceed the error term of (1.21).

This completes the proof of (1.21) when $u \in \mathcal{D}_{J+\nu}(y)$ and $\kappa = \nu \in \mathbb{N}^*$.

When $u \notin \mathcal{D}_{J+\nu}(y)$, $u > 1$, we have $\ell < u \leqslant \ell + 1$ for some integer $\ell \in [1, J+\nu+1]$. The Taylor-Lagrange formula must then take the contributions of the discontinuities into account. The quantity

$$\mathfrak{R}_1 := \sum_{1 \leqslant j \leqslant J+\kappa+1} \frac{(-1)^{j+1}}{j!} \sum_{\substack{1 \leqslant m \leqslant j \\ u-v < m < u}} \delta_{\kappa,m,j}(v+m-u)^j$$

is to be added to $\mathfrak{R}_0$. This is proved identically to [4; lemma 4.2] and we omit the details.[4] Arguing as before, we see that the quantity

$$\mathfrak{S}_1 := \sum_{1 \leqslant j \leqslant J+\kappa+1} \frac{(-1)^{j+1}}{j!} \sum_{\substack{1 \leqslant m \leqslant j \\ u-\varepsilon_y/2 < m < u}} \delta_{\kappa,m,j} \int_{u-m}^{\varepsilon_y/2} (v+m-u)^j \,\mathrm{d}\mu_{f,y}(v)$$

should be added to $\mathfrak{S}_0$. By (1.26), we see that the contribution of $m \neq \ell$ to $\mathfrak{S}_1$ is

$$\ll \frac{x}{L_{\beta+\delta/2}(x)}$$

and therefore may be absorbed by the error term of (1.21). When $m = \ell \geqslant 1$, we may again extend the integrals in $\mathfrak{S}_1$ to infinity and so obtain the term $U_J(x, y; f)$ appearing in (1.21).

It remains to evaluate $I_2$. That the order of magnitude of this quantity does not exceed that of the remainder of (1.21) readily follows from (5.5) via partial summation on observing that $\varphi_\kappa^{(j)}(v) = 0$ for $v \leqslant 1$ and $j \geqslant 1$ when $\kappa \in \mathbb{N}^*$.

### 5.2. *The case $\kappa \notin \mathbb{N}^*$*

With previous notation, we now have $\kappa = \nu + \vartheta$, $\nu \in \mathbb{N}$, $\vartheta \in\, ]0,1[$.

Our proof requires to a number of lemmas which we now state. Writing $v^+ := \max(0, v)$, we define

$$g_\vartheta(v) := \frac{1}{2\pi i} \int_{1+i\mathbb{R}} \frac{\mathrm{e}^{vs}}{s^\vartheta(s+1)}\,\mathrm{d}s = \frac{1}{\Gamma(\vartheta)} \int_0^{v^+} \mathrm{e}^{t-v} t^{\vartheta-1}\,\mathrm{d}t \qquad (v \in \mathbb{R}).$$

---

4. Note however that since we assumed $u \geqslant 1$ and $0 \leqslant v \leqslant \tfrac{1}{2}$, only the discontinuities at integers $m \geqslant 1$ appear in $\mathfrak{R}_1$.

**Lemma 5.1.** *For $\sigma > 0$, $T \geqslant 1$, $\vartheta \in ]0,1[$, $v \in \mathbb{R}^*$, we have*

$$\frac{1}{2\pi i} \int_{\sigma-iT}^{\sigma+iT} e^{sv} \frac{ds}{s^\vartheta} = \frac{v^{\vartheta-1}e^{\sigma v}}{\Gamma(\vartheta)} + O\Big(\frac{v^{\vartheta-1}e^{\sigma v}}{(1+|v|T)^\vartheta}\Big), \tag{5.6}$$

$$\frac{1}{2\pi i} \int_{\sigma-iT}^{\sigma+iT} \frac{e^{vs}}{s^\vartheta(1+s)} \, ds = g_\vartheta(v) + O\Big(\frac{e^{\sigma v}}{T^\vartheta(1+|v|T)}\Big). \tag{5.7}$$

*where the implicit constants depend at most on $\vartheta$.*

*Proof.* Estimate (5·6) follows from the bound

$$\int_T^\infty e^{sv} \frac{d\tau}{s^\vartheta} \ll \frac{|v|^{\vartheta-1}e^{\sigma v}}{(1+|v|T)^\vartheta} \quad (v \in \mathbb{R}^*). \tag{5.8}$$

When $|v| > 1/T$, partial integration yields (5·8). When $T \leqslant 1/|v|$ we may write

$$\int_T^{1/|v|} e^{i\tau v} \frac{d\tau}{s^\vartheta} = \int_T^{1/|v|} (e^{i\tau v} - 1)\frac{d\tau}{s^\vartheta} + \int_T^{1/|v|} \frac{d\tau}{s^\vartheta} \ll |v|^{\vartheta-1},$$

where we used the fact that $|e^{i\tau v} - 1| \leqslant |\tau v|$.

Estimate (5·7) is proved in [6; lemma 3.3]. □

For the following statement, we define $\|t\| := \min_{n \in \mathbb{Z}} |t-n|$ $(t \in \mathbb{R})$.

**Lemma 5.2.** *Let $\beta > 0$, $\mathfrak{c} > 0$, $\delta > 0$, $\kappa_0 > 0$, $\kappa > 0$, $\beta + \delta < 3/5$, $f \in \mathcal{H}(\kappa, \kappa_0; \beta, \mathfrak{c}, \delta)$. The integral*

$$Z_{f,y}(v) := \frac{(\log y)^{1-\vartheta}}{2\pi i} \int_{1+i\mathbb{R}} \mathcal{F}(1+s) y^{vs} \frac{s^{1-\vartheta}}{s+1} \, ds \tag{5.9}$$

*converges almost everywhere, the exceptional set being included in $\{v \in \mathbb{R} : y^v \in \mathbb{N}\}$. Moreover, for any fixed $j \geqslant 0$,*

$$\int_v^\infty t^j Z_{f,y}(t) \, dt \ll \frac{1}{(\log y)^{\vartheta+j} L_{\beta+\delta/2}(y^v)^2}, \tag{5.10}$$

*and, when $\|y^v\| \gg 1$, $0 \leqslant h \leqslant \frac{1}{2}$,*

$$\int_{v-h}^v Z_{f,y}(t) \, dt \ll h^\vartheta L_{\beta+\delta/2}(y^v/h \log y)^{2\vartheta-2}. \tag{5.11}$$

*Proof.* The statement regarding convergence is implied by (5·6). Thus

$$\int_0^v Z_{f,y}(w) \, dw = \frac{1}{2\pi i (\log y)^\vartheta} \int_{1+i\mathbb{R}} \mathcal{F}(1+s) \frac{y^{vs}-1}{s^\vartheta(s+1)} \, ds. \tag{5.12}$$

Splitting the integrand in (5·12) by isolating the term involving $y^{vs}$ and moving the integration line for this part into a zero-free region of $\zeta(1+s)$, we get

$$\int_0^v Z_{f,y}(t) \, dt = \frac{-1}{2\pi i (\log y)^\vartheta} \int_{1+i\mathbb{R}} \frac{\mathcal{F}(1+s)}{s^\vartheta(s+1)} ds + O\Big(\frac{(\log y)^{-\vartheta}}{L_{\beta+\delta/2}(y^v)^2}\Big),$$

which implies (5·10) for $j = 0$. The extension to $j \geqslant 1$ is immediate.

The proof of (5·11) is more delicate and follows the approach of [6; lemma 3.4]. We observe at the outset that we may assume $h \leqslant 1/L_{\beta+\delta/2}(y^v)$ since (5·11) otherwise follows from (5·10). To simplify the exposition, we prove (5·11) in the case $y := e$, and consequently assume $h \leqslant 1/L_{\beta+\delta/2}(e^v)$. We start with a general upper bound for

$$b_f(v) := \int_{1+i\mathbb{R}} \frac{\mathcal{F}(1+s)e^{vs}}{s^\vartheta(s+1)} \, ds$$

in which $\|e^v\|$ need not be bounded from below. Let $N \in \mathbb{N}^*$ be such that $\|e^v\| = |e^v - N|$. For $v \geqslant 3$, $n \in \mathbb{N}^*$, let us specialize in (5·7) $\sigma := 1/v$, replace $v$ by $v - \log n$, multiply out by $f(n)/n$ and sum up over $n \in \mathbb{N}^*$. We get

$$b_f(v) = \int_{1/v-iT}^{1/v+iT} \frac{\mathcal{F}(1+s)e^{vs}}{s^\vartheta(s+1)}\,\mathrm{d}s + O\bigg(\frac{1}{T^\vartheta}\sum_{n\geqslant 1}\frac{f^\dagger(n)}{n^{1+1/v}(1+T|v-\log n|)}\bigg).$$

By (4·6), we see that the contribution to the last sum of those integers $n$ outside the interval $[N/2, 3N/2]$ is $\ll v^{\kappa_0}/T$. That of the term $n = N$ is, still by (4·7),

$$\ll \frac{e^v}{(e^v + T\|e^v\|)L_{\beta+\delta/2}(e^v)^{10}}.$$

Let us denote by $V$ the complementary contribution, so that

$$V \ll \sum_{0<|n-N|\leqslant N/2} \frac{f^\dagger(n)}{n+T|N-n|}.$$

If $T > N$, this implies

$$V \ll \frac{1}{T} \sum_{0\leqslant m\leqslant L_{\beta+\delta/2}(N)^{10}} \frac{1}{m+1} \sum_{mN\leqslant |n-N|L_{\beta+\delta/2}(N)^{10}\leqslant (m+1)N} f^\dagger(n)$$

$$\ll \frac{N(\log N)^{\beta+\kappa_0}}{TL_{\beta+\delta/2}(N)^{10}} \ll \frac{1}{L_{\beta+\delta/2}(e^vT)^4}.$$

If $T \leqslant N$, we have similarly

$$V \ll \sum_{0\leqslant m\leqslant T/2} \frac{1}{(m+1)N} \sum_{mN/T<|n-N|\leqslant (N+1)N/T} f^\dagger(n)$$

$$\ll \sum_{0\leqslant m\leqslant T/2} \frac{1}{(m+1)N}\bigg\{\frac{N(\log N)^{\kappa_0-1}}{T} + \frac{N}{L_{\beta+\delta/2}(N)^{10}}\bigg\} \ll \frac{v^{\kappa_0}}{T} + \frac{1}{L_{\beta+\delta/2}(e^vT)^4}.$$

Hence, in all circumstances, writing $\beta_1 := \beta + \delta/2$,

$$b_f(v) = \int_{1/v-iT}^{1/v+iT} \frac{\mathcal{F}(1+s)e^{vs}}{s^\vartheta(s+1)}\,\mathrm{d}s + O\bigg(\frac{v^{\kappa_0}}{T^{1+\vartheta}} + \frac{1}{T^\vartheta L_{\beta_1}(e^vT)^4} + \frac{e^v}{T^\vartheta(e^v+T\|e^v\|)L_{\beta_1}(e^v)^{10}}\bigg).$$

The same estimate holds for $b_f(v-h)$ at the cost of replacing $v$ by $v-h$ in the integrand and $\|e^v\|$ by $\|e^{v-h}\|$ in the remainder. Therefore, if $\|e^v\| \gg 1$,

$$b_f(v) - b_f(v-h) \ll h\int_{1/v-iT}^{1/v+iT} \frac{|\mathcal{F}(1+s)||s|^{1-\vartheta}}{|s+1|}|\,\mathrm{d}s| + \frac{v^{\kappa_0}}{T^{1+\vartheta}}$$

$$+ \frac{1}{T^\vartheta L_{\beta_1}(e^vT)^4} + \frac{e^v}{T^\vartheta(e^v+T\|e^{v-h}\|)L_{\beta_1}(e^v)^{10}}$$

$$\ll hT^{1-\vartheta}(\log T)^\kappa + \frac{v^{\kappa_0}}{T^{1+\vartheta}} + \frac{1}{T^\vartheta L_{\beta_1}(e^vT)^4}$$

$$+ \frac{e^v}{T^\vartheta(e^v+T\|e^{v-h}\|)L_{\beta_1}(e^v)^{10}}.$$

Now, there exists a constant $c > 0$ such that $\|e^{v-h}\| \gg 1$ if $h \leqslant ce^{-v}$. In this case, the last term of the above upper bound is

$$\ll \frac{e^v}{T^{1+\vartheta}L_{\beta_1}(e^v)^{10}} \ll \frac{1/h}{T^{1+\vartheta}L_{\beta_1}(e^v/h)^5}.$$

If, to the contrary, $h > ce^{-v}$, then this last term is

$$\ll \frac{1}{T^\vartheta L_{\beta_1}(e^v)^{10}} \ll \frac{1}{T^\vartheta L_{\beta_1}(e^v/h)^5}.$$

Thus, in all cases,

$$b_f(b) - b_f(v-h) \ll hT^{1-\vartheta}(\log T)^\kappa + \frac{v^{\kappa_0}}{T^{1+\vartheta}} + \frac{1}{T^\vartheta L_{\beta+\delta/2}(e^vT)^4} + \frac{1+1/(hT)}{T^\vartheta L_{\beta+\delta/2}(e^v/h)^4}.$$

Selecting $T := 1/\{hL_{\beta+\delta/2}(e^v/h)^2\}$, we get (5·11) for $y = e$ and hence for general $y$ on replacing $v$ by $v\log y$ and $h$ by $h\log y$. □

**Lemma 5.3.** *Under the hypotheses of Lemma 5.2 and with the convention that $a_h(f) = 0$ if $h < 0$, we have*

$$\frac{(-1)^j}{j!} \int_0^\infty Z_{f,y}(v) v^j \, dv = \frac{a_{j-\nu-1}(f)}{(\log y)^{j+\vartheta}} \quad (j \geqslant 0). \tag{5.13}$$

*Proof.* The bound (5·10) guarantees via partial integration that the integrals (5·13) converge. We have

$$\widehat{Z_{f,y}}(s) = \frac{s^{1-\vartheta} \mathcal{F}(1 + s/\log y)}{s + \log y}, \tag{5.14}$$

and observe that this implies, for $\vartheta \neq 0$,

$$d\mu_{f,y}(v) = Z_{f,y}(v) \, dv. \tag{5.15}$$

From (5·14), we get, for $|s| < \mathfrak{c} \log y$, $J \geqslant 1$,

$$\widehat{Z_{f,y}}(s) = \frac{s^{\nu+1} \mathcal{F}(1 + s/\log y)}{s^\kappa (1 + s/\log y) \log y} = \sum_{0 \leqslant j \leqslant J} \frac{a_j(f) s^{j+\nu+1}}{(\log y)^{\kappa+j+1}} + O_J\left(\frac{s^{J+2+\nu}}{(\log y)^{\kappa+J+2}}\right)$$

$$= \int_0^\infty Z_{f,y}(v) e^{-vs} \, dv = \sum_{0 \leqslant j \leqslant J+\nu+1} \frac{(-1)^j s^j}{j!} \int_0^\infty Z_{f,y}(v) v^j \, dv + O_J\left(\frac{s^{J+2+\nu}}{(\log y)^{\kappa+J+2}}\right).$$

This is all we need. □

We are now in a position to complete the proof of (1·21) when $\kappa \in \mathbb{R}^+ \smallsetminus \mathbb{N}$. In view of (4·7), we may plainly assume that $x = y^u \in \frac{1}{2} + \mathbb{N}^*$.

Since, from (5·14) and (3·5),

$$\frac{\mathcal{B}(s)}{s\zeta(s)^\kappa s_y^\kappa \widehat{\varrho}(s_y)^\kappa} = \frac{s_y \mathcal{B}(s)}{s\zeta(s)^\kappa s_y^\vartheta} \frac{1}{s_y^{\nu+1} \widehat{\varrho}(s_y)^\kappa} = \widehat{Z_{f,y}}(s_y) \widehat{\varphi_\kappa}(s_y) \log y.$$

we deduce by (4·12) that

$$A(x, y; f) = J_y(\log x) = x(Z_{f,y} * \varphi_\kappa)(u). \tag{5.16}$$

The above convolution may be estimated in much the same way as in the case $\kappa \in \mathbb{N}^*$. We write

$$\int_0^u \varphi_\kappa(u-v) Z_{f,y}(v) \, dv = I_1 + I_2 + I_3, \tag{5.17}$$

where $I_1$ corresponds to the contribution of $v \in [0, \frac{1}{2}\varepsilon_y]$, $I_2$ to that of $v \in ]\frac{1}{2}\varepsilon_y, u - \frac{1}{2}]$, and $I_3$ to that of $[u - \frac{1}{2}, u]$.

Partial integration yields $I_2 \ll 1/L_{\beta+\delta/2}(y^{\varepsilon_y/2})$ in view of (5·10). To bound $I_3$, we use (5·11) and (1·11). Consequently

$$I_3 \ll \left[(u-v)^{-\vartheta} \left| \int_v^u Z_{f,y}(t) \, dt \right| \right]_{u-1/2}^u + \int_0^{1/2} \frac{dh}{hL_{\beta+\delta/2}(y^u/h \log y)^{2-2\vartheta}} \ll \frac{1}{L_{\beta+\delta/2}(y^u)^{1-\vartheta}}.$$

To estimate $I_1$, we argue differently according to whether $u \in \mathcal{D}_{J+\nu}(y)$ or not. In the first instance, the Taylor expansion (5·3) is still valid and we derive the required conclusion using (5·10), (5·13) and (5·5).

In the complementary case, let $\ell \in [1, J + \nu + 1]$ be defined by $\ell < u \leqslant \ell + 1$. Writing the Taylor expansion at order $\ell$ provides

$$\varphi_\kappa(u-v) = \sum_{0 \leqslant j < \ell} \frac{(-v)^j \varphi_\kappa^{(j)}(u)}{j!} + \frac{(-1)^\ell}{(\ell-1)!} \int_0^v (v-t)^{\ell-1} \varphi_\kappa^{(\ell)}(u-t) \, dt.$$

The contribution of the sum to (5·16) is well approximated by

$$\sum_{0 \leqslant j < \ell - \nu - 1} \frac{a_j(f)\psi_\kappa^{(j+1)}(u)}{(\log y)^{\kappa+j+1}},$$

this sum being empty if $\ell \leqslant \nu + 1$, that is $u \leqslant \nu + 2$. The contribution of the integral to $I_1$ is

$$\frac{(-1)^\ell}{(\ell-1)!} \int_0^{\varepsilon_y/2} \varphi_\kappa^{(\ell)}(u-t) \int_t^{\varepsilon_y/2} (v-t)^{\ell-1} Z_{f,y}(v) \, \mathrm{d}v \, \mathrm{d}t,$$

which, to an admissible error, is approximated by $X_\ell(x, y; f)$, as defined in (1·18). Estimate (1·20) is obtained from (5·10) by bounding $\varphi_\kappa^{(\ell)}(u-t)$ by $\ll (u - \ell - t)^{-\vartheta}$ for $0 \leqslant t \leqslant u - \ell$, and by $\ll 1$ for $u - \ell \leqslant t \leqslant \tfrac{1}{2}\varepsilon_y$.

This finishes the proof.

## 6. Proof of Theorem 2.1

We have, for positive, non-integral $x$,

(6·1) $$m(x, y; f) = \frac{1}{2\pi i(\log y)^\kappa} \int_{1+i\mathbb{R}} \frac{\mathcal{B}(s+1)x^s}{\zeta(s+1)^\kappa s^{\kappa+1} \widehat{\varrho}(s \log y)^\kappa} \, \mathrm{d}s.$$

Arguing as in the proof of Theorem 1.1, we see that, to the stated error, a suitable approximation to $m(x, y; f)$ is

$$\int_\mathbb{R} h_\kappa(u-v) y^{-v} \, \mathrm{d}M(y^v; f).$$

Writing $m(x; f) := m(x, x; f)$, we have $y^{-v} \, \mathrm{d}M(y^v; f) = \mathrm{d}m(y^v; f)$.

When $\kappa \in \mathbb{N}^*$, we may integrate the Taylor-Maclaurin expansion (5·3) taking the contribution of the discontinuities into account and get (2·6).

When $\kappa \in \mathbb{R}^+ \smallsetminus \mathbb{N}$, we rewrite (6·1) as

$$m(x, y; f) = \frac{(\log y)^{1-\vartheta}}{2\pi i} \int_{1+i\mathbb{R}} \frac{\mathcal{B}(s+1)\widehat{\varphi_\kappa}(s \log y) x^s}{\zeta(s+1)^\kappa s^\vartheta} \, \mathrm{d}s,$$

whence, for suitable $\varepsilon > 0$, as in (4·10),

$$m(x, y; f) = (Z_{f,y}^* * \varphi_\kappa)(u) + O\left(\frac{R_\kappa(u)}{L_\varepsilon(y)}\right)$$

with

$$Z_{f,y}^*(v) := \frac{1}{2\pi i} \int_{1+i\mathbb{R}} \frac{\mathcal{F}(1 + s/\log y) \mathrm{e}^{vs}}{s^\vartheta} \, \mathrm{d}s.$$

The estimates of Lemma 5.2 remain valid for $Z_{f,y}^*$ and we have, parallel to (5·13),

$$\frac{(-1)^j}{j!} \int_\mathbb{R} v^j Z_{f,y}^*(v) \, \mathrm{d}v = \frac{a_{j-\nu}^*(f)}{(\log y)^{j+\vartheta}},$$

with the same convention as in Lemma 5.3 regarding negative values of the index.

Next, recalling notation (1·12), we apply (5·3) and restrict the summation to $j \leqslant J + \nu$. Handling as before the contribution of the range $\tfrac{1}{2} \leqslant v \leqslant u$ to the integral $(Z_{f,y}^* * \varphi_\kappa)(u)$ as an error term, we get the desired estimates in Theorem 2.1.

## 7. Proof of Corollary 1.4

We first observe that [15; th. 2.1] provides in $G_\beta$, with any $\beta < 1$,

(7·1) $$M(x,y;f^\dagger) \asymp \mathcal{B}(\alpha_\kappa,y)^{O(1)} x \varrho_\kappa(u)(\log y)^{\kappa-1}.$$

The proof of Corollary 1.4 closely follows that of [12; th. 1], corresponding to the case $f = \mu$, $f^\dagger = \mathbf{1}$. We distinguish several cases according to the size of $u$.

When $1 \leqslant u \leqslant \nu + 2$, the required estimate readily follows from (1·24) and (7·1).

When $u > \nu + 2$, $\beta < 3/5$, $(x,y) \in G_\beta$, we appeal to (1·9), (1·21) and (3·13) to get the upper bound

$$\ll \frac{e^{-c_1 u/(\log 2u)^2}}{(\log y)^{2\kappa}},$$

valid for any $c_1 < \pi^2/2$.

When $u > (\log y)^{\beta/(1-\beta)}$, we redefine $\alpha_\kappa$ as the saddle point associated to the Perron integral for $M(x,y;\tau_\kappa)$, hence involving $\zeta(s,y)^\kappa$. We then have, by Perron's formula, arguing as in [12; lemma 2], for $r < 3/2$,

$$M(x,y;f) = \frac{1}{2\pi i} \int_{\alpha_\kappa - iL_r(y)^2}^{\alpha_\kappa + iL_r(y)^2} \frac{\mathcal{B}(s,y)x^s}{s\zeta(s,y)^\kappa}\,ds + O\Big(\frac{x\varrho_\kappa(u)}{L_r(y)} + x\varrho_\kappa(u)e^{-c_2 u}\Big).$$

We conclude following the proof of [12; th. 1] by appealing to the bound [12; (2.4)] for $\zeta(s,y)/\zeta(\alpha_\kappa,y)$ and to the saddle point estimate stated in [3; §2.1] for $M(x,y;\tau_\kappa)$, which extends to $M(x,y,f^\dagger)$.

## 8. Proof of Theorem 2.2

For $f$ as in the statement, we put $g = f * \mu$, so that $g \in \mathcal{H}(\kappa+1,\kappa_0+1;\beta,\mathfrak{c},\delta)$. We may plainly assume $y$, and hence $x$, sufficiently large throughout the proof.

Let $z$ be a parameter to be defined later. By the hyperbola principle, we have

$$M(x;f_y) = \sum_{n \in S(x/z,y)} g(n)\Big\lfloor\frac{x}{n}\Big\rfloor + \sum_{d \leqslant z} M\Big(\frac{x}{d},y;g\Big) - M\Big(\frac{x}{z},y;g\Big)\lfloor z \rfloor =: S_1 + S_2 - S_3.$$

We can replace $\lfloor x/n \rfloor$ by $x/n$ in $S_1$ with an error not exceeding that of (2·8) granted that

$$\frac{x\varrho_{\kappa+1}(u)(\log y)^\kappa}{z} \ll \frac{xR_{\kappa+1}(u)}{(\log y)^{J+1}}.$$

This is certainly the case for the choice $z := e^u(\log y)^{\kappa+J+2}$. Under the assumptions of the statement, we then have $z = y^{o(1)}$ and we may apply Theorem 2.1 to $g$, getting the suitable approximation

$$S_1 \approx x \sum_{0 \leqslant j \leqslant J} a_j^*(g)\frac{\psi_{\kappa+1}^{(j)}(u - u_z)}{(\log y)^{\kappa+j+1}}.$$

Here and in the sequel of this proof we use the symbol $A \approx B$ to indicate that two quantities $A$ and $B$ agree to within an error not exceeding that of (2·8).

By Theorem 1.1 and (5·16), we have, writing $u_t := (\log t)/\log y$,

(8·1) $$S_2 \approx \int_0^{1/2} \sum_{d \leqslant z} \frac{\varphi_{\kappa+1}(u - v - u_d)}{d}\,d\mu_{g,y}(v).$$

The inner sum may be rewritten as

$$\int_1^z \varphi_{\kappa+1}(u - v - u_t)\frac{d\lfloor t \rfloor}{t} = V_{21} + V_{22},$$

with

$$V_{21} := (\log y)\int_0^{u_z} \varphi_{\kappa+1}(u - v - w)\,dw, \qquad V_{22} := -\int_1^z \varphi_{\kappa+1}(u - v - u_t)\frac{d\langle t \rangle}{t}.$$

The contribution of $V_{21}$ to (8·1) is

$$\approx (\log y)\sum_{0 \leqslant j \leqslant J + \nu + 2}\frac{(-1)^j}{j!}\int_0^{u_z}\varphi_{\kappa+1}^{(j)}(u - w)\,dw \int_0^\infty v^j\,d\mu_{g,y}(v)$$

$$= \sum_{0 \leqslant j \leqslant J}\frac{a_j(g)\{\psi_{\kappa+1}^{(j)}(u) - \psi_{\kappa+1}^{(j)}(u - u_z)\}}{(\log y)^{\kappa+j+1}},$$

where Lemma 5.3 and (5·11) have been used.

The quantity $V_{22}$ is handled by writing

$$V_{22} \approx \sum_{0 \leqslant h \leqslant J+\nu+2} \frac{(-1)^{h+1}}{h!} \varphi_{\kappa+1}^{(h)}(u-v) \int_1^z u_t^h \frac{\mathrm{d}\langle t \rangle}{t} \approx \sum_{0 \leqslant h \leqslant J+\nu+2} \frac{\gamma_{h+1} \varphi_{\kappa+1}^{(h)}(u-v)}{(\log y)^h},$$

where $\{\gamma_h\}_{h=0}^{\infty}$ is the sequence of Taylor coefficients of $s\zeta(1+s)$ at the origin. Carrying back into (8·1), we obtain the contribution

$$\sum_{0 \leqslant h \leqslant J+\nu+2} \frac{\gamma_{h+1}}{(\log y)^h} \sum_{0 \leqslant j \leqslant J+\nu+2} \frac{(-1)^j \varphi_{\kappa+1}^{(h+j)}(u)}{j!} \int_0^{\infty} v^j \, \mathrm{d}\mu_{g,y}(v)$$

$$\approx \sum_{0 \leqslant h \leqslant J} \sum_{0 \leqslant j \leqslant J} \frac{\gamma_{h+1} a_j(g) \psi_{\kappa+1}^{(h+j+1)}(u)}{(\log y)^{h+j+\nu+2+\vartheta}}$$

$$\approx \sum_{0 \leqslant h \leqslant J} \sum_{0 \leqslant j \leqslant J} \frac{\gamma_{h+1} a_j(g) \psi_{\kappa+1}^{(h+j+1)}(u)}{(\log y)^{\kappa+2+h+j}} \approx \sum_{0 \leqslant j \leqslant J} \frac{\psi_{\kappa+1}^{(j+1)}(u)\{a_{j+1}(f) - a_{j+1}(g)\}}{(\log y)^{\kappa+2+j}}.$$

Finally, we have

$$S_3 \approx x \sum_{0 \leqslant j \leqslant J} \frac{a_j(g) \psi_{\kappa+1}^{(j+1)}(u-u_z)}{(\log y)^{\kappa+j+2}}.$$

Gathering our estimates and using the fact that $a_0(g) = a_0(f)$, we obtain (2·8) as required.

Régis de la Bretèche
Université Paris Cité, Sorbonne Université, CNRS
Institut de Math. de Jussieu-Paris Rive Gauche
F-75013 Paris
France

regis.de-la-breteche@imj-prg.fr

Gérald Tenenbaum
Institut Élie Cartan
Université de Lorraine
BP 70239
54506 Vandœuvre-lès-Nancy Cedex
France

gerald.tenenbaum@univ-lorraine.fr